\documentclass[a4paper, 12pt]{amsart}
\numberwithin{equation}{section}
\usepackage{amssymb}
\usepackage{amsmath}
\usepackage{amsfonts}
\usepackage{marvosym}
\usepackage{tikz}

\usepackage{stackrel}

\newtheorem{theorem}{Theorem}[section]
\newtheorem{thm}{Theorem}[section]
\newtheorem{lem}[theorem]{Lemma}

\newtheorem{prop}[theorem]{Proposition}

\theoremstyle{definition}

\newtheorem{rem}[theorem]{Remark}

\def\cb{{\mathcal B}}

\def\ga{{\mathfrak A}}
\def\gb{{\mathfrak B}}
\def\gc{{\mathfrak C}}

\def\gn{{\mathfrak N}}

\def\gx{{\mathfrak X}}

\def\bc{{\mathbb C}}
\def\bd{{\mathbb D}}

\def\bm{{\mathbb M}}
\def\bn{{\mathbb N}}

\def\br{{\mathbb R}}
\def\bt{{\mathbb T}}

\def\bz{{\mathbb Z}}

\def\b{\beta}
\def\g{\gamma}  
\def\d{\delta}

\def\l{\lambda} 

\def\m{\mu}

\def\s{\sigma} 

  \def\F{\Phi}

\def\id{\hbox{id}}

\def\dist{\mathop{\rm dist}}

\def\di{{\rm d}}
\def\id{{\rm id}}

\def\idd{{1}\!\!{\rm I}}

\DeclareMathAlphabet{\mathpzc}{OT1}{pzc}{m}{it}

\begin{document}

\title[]
{Decoherence for Markov chains}
\author{Francesco Fidaleo}
\address{Francesco Fidaleo\\
Dipartimento di Matematica \\
Universit\`{a} di Roma Tor Vergata\\
Via della Ricerca Scientifica 1, Roma 00133, Italy} \email{{\tt
fidaleo@mat.uniroma2.it}}
\author{Elia Vincenzi}
\address{Elia Vincenzi\\
Dipartimento di Matematica \\
Universit\`{a} di Roma Tor Vergata\\
Via della Ricerca Scientifica 1, Roma 00133, Italy} \email{{\tt
vincenzi@mat.uniroma2.it}}
\date{\today}

\begin{abstract}
The subspace generated by the eigenvectors pertaining to the peripheral spectrum of any stochastic matrix is naturally equipped with a structure of a (finite dimensional abelian) $C^*$-algebra, and the restriction of such a stochastic matrix to this subspace, indeed a $C^*$-algebra under this canonical new product, generates a conservative $C^*$-dynamical system. \end{abstract}

\maketitle

\section{introduction}
The universally accepted concept of {\it decoherence} was codified by P. Blanchard and R. Olkiewicz, and concerns the properties of separation in the persistent and transient parts of a dissipative $C^*$-dynamical system, typically a "small" system which interacts with a huge reservoir, see \cite{BO}. 

Such kind of $C^*$-dynamical systems are described by a $C_o$ semigroup but, in order to capture most of the main properties, we can also consider discrete dynamics generated  by linear Unital Completely Positive (UCP for short) maps. In all the forthcoming analysis, we restrict the matter to the simpler picture described by such maps.

Indeed, let $(\ga,\F)$ be a $C^*$ dynamical system consisting of a unital $C^*$-algebra $\ga$ on which the UCP map $\F$ is acting. In this context, it is said that the decoherence takes place if the whole system is split in some kind of direct sum as  $\ga=\gn(\F)\bigoplus \ga_o$,
where 
$$
\gn(\F):=\{x\in\ga\mid \F(x^*x)=\F(x^*)\F(x),\,\F(xx^*)=\F(x)\F(x^*)\}\,,
$$
often referred to as the multiplicative domain of $\F$,
and 
$$
\ga_o:=\Big\{x\in\ga\mid\lim_{n\to+\infty}\|\F^n x\|=0\Big\}\,.
$$

As it is explained in \cite{FOR}, the occurrence of any reasonable meaning of decoherence seems to be strictly connected to the spectral properties of the involved UCP map $\F$. For the majority of the cases of interest, but the case of $*$-automorphisms for which the decoherence is trivially satisfied, $\s(\F)$ is the whole disk, and thus it appears complicated to provide a splitting of $\ga$ in a direct sum as above.\footnote{If one deals with $C_o$ semigroup $e^{-tA}$, the case of the whole closed unit disc for the spectrum of $\F$, would correspond for the generator $A$ to $\s(A)$ being equal to the whole right half plane.}

%It is of certain interest to note that all UCP maps acting on finite dimensional $C^*$-algebras are gapped
However, for the so called {\it gapped} UCP maps, which are those where the peripheral spectrum is topologically separated from the part inside the unit disk, a splitting of the involved algebra 
$$
\ga=P_\F\ga+(\idd_{\cb(\ga)}-P_\F)\ga
$$ 
in the persistent and transient part is easily obtained using the holomorphic functional calculus, see  \cite{FOR}.

It is also straightforward to see that the persistent part $P_\F\ga$ is a norm-closed subspace, trivially containing the identity of $\ga$, in general only closed by taking the adjoint and not under the product, that is it is merely an {\it operator system}, see {\it e.g.} \cite{CE}.

Among some interesting situations of gapped  UCP maps, there are certainly all those acting on finite dimensional $C^*$-algebras. In this, apparently simple, situation for which a quite complete description of properties of UCP maps is however a formidable task, we want to point out that the persistent part can be always equipped with a new product ({\it cf.} \cite{CE}, Theorem 3.1), in general different from the original one, making such a persistent part a "genuine" $C^*$-algebra.

It is now to address the natural question whether the UCP map $\F$ under consideration, restricted to $P_\F\ga$, the latter equipped with this canonical new product, is indeed an automorphism and then providing a genuine conservative dynamical system. In this way, the decoherence takes place just by considering the persistent part being equipped with its natural structure of a $C^*$-algebra. In this regard, it should be noticed that some simple and natural examples do not enjoy the property of decoherence according to the commonly used definition described above, see \cite{FOR}, Section 6.

The solution to this interesting problem is still open, even in the simplest finite dimensional case. Yet, it would be of interest to provide an answer, at least for some specified class of UCP maps. Among those class of dissipative dynamical systems there are those associated to {\it Markov chains}, the last being among ones of the most interesting, and therefore most studied, objects in mathematics.

The present note is indeed devoted to prove that any Markov chain is encoding a genuine conservative dynamical system, after separating the persistent part from the transient part, exponentially disappearing in the limit taken on the iteration of the involved stochastic matrix or, concretely after a finite number of such interactions depending on the nature of the concrete model and the size of the so-called {\it mass gap}, which is the distance between the peripheral spectrum of the stochastic matrix, and the part of the spectrum inside the unit disk.

\section{preliminaries}

\medskip

\noindent
\textbf{Basic notation.} 
With $\bd$ we denote the closed unit disk $\{\l\in\bc\mid|\l|\leq1\}$. Its boundary $\partial\bd$ is nothing else than the unit circle denoted by $\bt$.

If it is not otherwise specified, in the present note we only deal with everywhere defined linear maps between vector spaces. All $C^*$-algebras $\ga$ we deal with are unital with 
$I=\idd_\ga$. For the Banach algebra $\cb(\ga)$ consisting of all bounded operators acting on the $C^*$-algebra $\ga$, we also put $I=\idd_{\cb(\ga)}$.

For involutive algebras $\gc_i$, $i=1,2$, a map $\Psi:\gc_1\to\gc_2$ is said to be {\it selfadjoint} or {\it real} if $\Psi(x^*)=\Psi(x)^*$ for every $x\in\gc_1$.

For the $C^*$-algebra $\ga$, the map $\F:\ga\to\ga$ is said {\it completely positive} if $\F\otimes\id_{\bm_n(\bc)}=:\F_n:\bm_n(\bc)\to\bm_n(\bc)$ is positive for each $n=1,2,\dots$\,. Thus,  it is positive if $\F_1$ is. It is {\it unital} if $\F(I)=I$. If $\ga$ is abelian, then complete positivity coincides with positivity, see {\it e.g.} \cite{St}, Theorem 1.2.4.

For a self-map $T\in\cb(\gx)$ on the Banach space $\gx$, the peripheral spectrum is $\s_{\rm ph}(T):=\{\l\in\s(T)\mid |\l|={\rm spr}(T)\}$, ${\rm spr}(T)$ being the spectral radius. Such a map $T$ is said to be {\it gapped} if it presents the so-called {\it mass-gap}, that is 
$\dist\big(\s_{\rm ph}(T),\s(T)\smallsetminus\s_{\rm ph}(T)\big)>0$.\footnote{The terminology "mass-gap" comes from physical motivations because such a distance is dimensionally equivalent to a mass.} 

In the case of gapped UCP maps $\F$, ${\rm spr}(\F)=1=\|\F\|$, and thus $\ga=P_\F\ga+Q_\F\ga$ with $Q_\F:=\frac1{2\pi\imath}\oint_\g\big(\l I-\F\big)^{-1}\di\l$, $P_\F:=I-Q_\F$, and $\g$ is a counterclockwise Jordan curve inside the open unit disk surrounding $\s(\F)\smallsetminus\s_{\rm ph}(\F)$.

Therefore, in the case of gapped UCP maps, $\lim_n\|\F^n(x)\|=0$, $x\in Q_\F\ga$,
see \cite{FOR}, Proposition 3.1.

\medskip

\noindent
\textbf{Jordan morphisms.} 
To reduce the matter to our setting, we deal only with Jordan algebras made of the selfadjoint part $\ga_{\rm sa}$ of $C^*$-algebras $\ga$ equipped with the Jordan product given by
$$
\ga_{\rm sa}\ni a,b\mapsto a\bullet b:=\frac12(ab+ba)\in\ga_{\rm sa}\,.
$$

A selfadjoint map $\F:\ga\to\gb$ between $C^*$-algebras $\ga$ and $\gb$ is said an {\it order-isomorphism} if it is invertible and $\F$ and $\F^{-1}$ are both positive. If $\ga$ coincides with $\gb$, we speak of {\it order-automorphism}.

For the convenience of the reader,  we report the following result, crucial for our analysis.
\begin{thm}[\cite{St}, Thm. 2.1.3]
\label{crstdh}
A unital selfadjoint map $\F:\ga\to\gb$ is a order-isomorphism if and only if $\F\lceil_{\ga_{\rm sa}}:\ga_{\rm sa}\to\gb_{\rm sa}$ is a Jordan isomorphism.
\end{thm}

\medskip

\noindent
\textbf{$C^*$-dynamical systems.}
A (discrete) $C^*$-dynamical system is simply a triple $(\ga,\F,M)$, where $\ga$ is a $C^*$-algebra, $\F$ a UCP map acting on $\ga$ in a natural way, and $M$ is the monoid $\bn$ or $\bz$. By definition, the case relative to the group $\bz$ corresponds to $\F$ being a $*$-automorphism, and in particular invertible. Therefore, this last case corresponds to (microscopically) reversible, $C^*$-dynamical systems.\footnote{The sentence "microscopically reversible" means that a system, typically describing temperature, or Kubo-Martin-Schwinger boundary condition, states, are notoriously irreversible in some sense, but microscopically reversible because their dynamics is generated by unitary operators, see e.g. \cite{BR}.} Reversible dynamical systems are often referred to as "conservative, Hamiltonian, or unitary".  The cases when $\F$ is not invertible, are referred to as {\it dissipative} dynamical systems.
Relatively to the generic case treated here when $\F$ is not a $*$-auto\-morphism, the simplified notation $(\ga,\F)$ stands for $(\ga,\F,\bn)$.

We end with the following
\begin{rem}
Let $\F:\ga\to\ga$ be an irreducible UCP map. By \cite{G}, Proposition 3.2, the closed subspace ${\rm AP}(\ga)$ of $\ga$ generated by the eigenvectors pertaining to the peripheral eigenvalues of $\F$ is a $C^*$-subalgebra of $\ga$, and the restriction $\F\lceil_{{\rm AP}(\ga)}$ of $\F$ to ${\rm AP}(\ga)$ is a $*$-automorphism.
\end{rem}
Here, ${\rm AP}(\ga)$ is known as the {\it almost periodic} part of $\F$, see {\it e.g.} \cite{F24}. In addition,
we are using the notion of irreducibility in \cite{G}, Definition 2.2. We see in Proposition \ref{vig} below that such a definition of irreducibility is equivalent to the analogous one used for positive matrices.

\medskip

\noindent
\textbf{Stochastic matrices.} 
A stochastic matrix $S\in\bm_n(\br)$ is nothing but a UCP map acting on $\bc^n$, and thus for all entries, $S_{ij}\geq0$. The $C^*$-dynamical systems
$(\bc^n, S)$ describe the dynamics of the so-called {\it Markov chains} ({\it e.g.} \cite{S}). The structure of a stochastic matrix is briefly outlined in Section \ref{stoc}. Here, we recall their basic properties.

Since the probability at each step of the transition ({\it i.e.} after the repeated application of $S$ on vectors of $\bc^n$) must be conserved, we have that 
$1:=\begin{bmatrix} 
	 1 \\
	1\\
	\cdot\\
	\cdot\\
	\cdot\\
	1\\
	\end{bmatrix}$
is a right eigenvector of $S$ with eigenvalue $1$: $S1=1$.

On the other hand, any positive row-vector $\begin{bmatrix} 
	 \pi_1 & \pi_2 &\cdot &\cdot&\cdot&\pi_n\\
	 \end{bmatrix}$
with $\sum_{i=1}^n\pi_i=1$, which is also a left eigenvector corresponding to the eigenvalue $1$, is a {\it stationary distribution} for the associated Markov chain. The multiplicity of the left and right eigenvalue 1 coincides.

A positive, and in particular a stochastic, matrix $A$ of order $n$ is said to be {\it irreducible} if there exists no permutation-matrix $P$ of order $n$ such that
$PAP^{-1}=\begin{bmatrix} 
	 A_1 & B\\
	O& A_2\\
	\end{bmatrix}$, where $A_i$ is a square-matrix of order $n_i$ with $1\leq n_i<n$, see {\it e.g.} \cite{S0, S}.
	
For the convenience of the reader, we show that the definition of irreducibility provided in Definition 2.2 of \cite{G}, coincides with the previous one for stochastic matrices.
\begin{prop}
\label{vig}
Let $A$ be a positive matrix. Then it is irreducible if and only if the only faces of the positive cone $\mathbb{C}^n_+$ which are invariant under $A$ are $\{0\}$ and the whole $\mathbb{C}^n_+$.
\end{prop}
\begin{proof}
Suppose that $A$ is irreducible, and there exists a face $F:=\bigoplus\limits_{j\in J}\mathbb{R}^+e_j$ for some $J\neq\varnothing,\{1,\dots,n\}$ which is invariant under $A$, that is $A(F)\subset F$. Consider any permutation $P$ which reorders the canonical basis as $\{e_{j}\}_{j\in J}\bigcup\{e_i\}_{i\in J^{\rm c}}$. Then,
$PAP^{-1}=\begin{bmatrix}
A_1&B\\
O&A_2
\end{bmatrix}$
with $A_1$ a square matrix of size $1\leq m_1=|J|\leq n-1$ and $A_2$ a square matrix of size $1\leq m_2\leq n-1$, which is a contradiction.

Suppose now that there exists a permutation matrix $P$ such that $PAP^{-1}=\begin{bmatrix}
A_1&B\\
O&A_2
\end{bmatrix}$ as above.
Consider the set $\{e_j\}_{j\in J}:=\{Pe_i\}_{i=1}^{m_1}$. Then,
$\bigoplus\limits_{j\in J}\mathbb{R}^+e_j$
is a non trivial face of $\mathbb{C}^n_+$ which is invariant under $A$, which is again a contradiction.
\end{proof}

\section{on UCP maps on finite dimensional $C^*$-algebras}

The present section is devoted to basic results, perhaps well known to the experts, on which is based the forthcoming analysis for stochastic matrices.
\begin{prop}
\label{prtreu}
Let $\F:\ga\to\ga$ be a UCP map on the finite dimensional $C^*$-algebra $\ga$. Then there exists a subsequence of natural numbers $(n_j)_j\subset\bn$ such that
$\lim_j \F^{n_j}=P_\F$.
\end{prop}
\begin{proof}
By using the Jordan decomposition of $\F$, we have
$$
\F=\sum_{\l\in\s_{\rm ph}(\F)}\l E_\l+Q_\F\F\,,
$$
with $E_\l E_\m=\d_{\l,\m}E_\l$, $\sum_{\l\in\s_{\rm ph}(\F)}E_\l=P_\F$. 

Since $\s_{\rm ph}(\F)$ is a finite subset of the circle (which not necessarily generates a finite subgroup), there exists a subsequence $(n_j)_j$ of natural numbers such that
$\lim_j\l^{n_j}=1$ for each $\l\in\s_{\rm ph}(\F)$. Then, by \cite{FOR}, Proposition 3.1, and the choice of the subsequence $(n_j)_j$, we get
\begin{align*}
\lim_j \F^{n_j}=&\lim_j\bigg(\sum_{\l\in\s_{\rm ph}(\F)}\l E_\l+Q_\F\F\bigg)^{n_j}\\
=&\lim_j\bigg(\sum_{\l\in\s_{\rm ph}(\F)}\l E_\l\bigg)^{n_j}+\lim_j\big(Q_\F\F^{n_j}\big)\\
=&\lim_j\sum_{\l\in\s_{\rm ph}(\F)}\l^{n_j}E_\l=\sum_{\l\in\s_{\rm ph}(\F)}\lim_j(\l^{n_j})E_\l\\
=&\sum_{\l\in\s_{\rm ph}(\F)}E_\l=P_\F\,.
\end{align*}
\end{proof}
We also report the well-known fact of the mean ergodicity of such maps.
\begin{prop}
Let $\F:\ga\to\ga$ be a UCP map on the finite dimensional $C^*$-algebra $\ga$. Then
$\lim_n\big(\frac1{n}\sum_{k=0}^{n-1}\F^k\big)=E_1$ the, necessarily UCP ({\it cf.} \cite{FM3}, Theorem 2.1), projection onto the fixed-point subspace $\ga^\F$.
\end{prop}
\begin{proof}
By performing the same calculations in the above proof, we get
$$
\lim_n\Big(\frac1{n}\sum_{k=0}^{n-1}\F^k\Big)=E_1+\sum_{\l\in\s_{\rm ph}(\F)\smallsetminus\{1\}}\frac1{1-\l}\lim_n\Big(\frac{1-\l^n}{n}\Big)E_\l
=E_1\,,
$$
since $|1-\l^n|\leq2$.
\end{proof}

Here, we report the main theorem concerning the injective operator systems relative to the finite dimensional situation.

For a UCP self-map $\F:\ga\to\ga$ of a finite dimensional $C^*$-algebra $\ga$, define $\ga_\F:=\big(P_\F\ga,\idd_\ga,*,\circ,\|\,\,\|_\ga\big)$ the operator system $P_\F\ga$ equipped with the binary operation
$a\circ b:=P_\F(ab)$.
\begin{thm}[Choi-Effros]
\label{cief}
The operator system $\ga_\F$, equipped with the new binary operation ''$\circ$'', is a $C^*$-algebra.
\end{thm}
\begin{proof}
See Theorem 3.1 of \cite{CE}.
\end{proof}

\section{stochastic matrices and the associated persistent dynamical systems}
\label{stoc}

The present section starts with useful facts on stochastic matrices. The second half is devoted to the main result of the present note.

After a possible permutation of indices, any such a matrix $S$ has the following canonical form (e.g.\ \cite{S, V})
\begin{equation}
\label{can}
S=\begin{bmatrix} 
	 B_{00} & B_{01} &\cdot &\cdot&\cdot&B_{0n}\\
	0& B_{11}& 0&\cdot&\cdot& 0\\
	\cdot& 0& B_{22}&0&\cdot& \cdot\\
	\cdot& \cdot& \cdot&\cdot&\cdot&\cdot\\
	\cdot& \cdot&\cdot&\cdot&\cdot&\cdot\\
	0&\cdot& \cdot&\cdot&\cdot& B_{nn}\\
	\end{bmatrix}\,. 
\end{equation}

Here,  $B_{00}$ is a square strictly sub-stochastic matrix associated to the transient indices, which can be the empty matrix if and only if the subset of such transient  indices is empty.\footnote{In the language of Markov chains, the transient indices are associated to the so called transient (or inessential) "states".} The  
square-block irreducible matrices $\{B_{11},\dots,B_{nn}\}$ describe the single ergodic components relative to the recurrent indices.

The above form for a stochastic matrix is said to be the {\it reduced form} and, if there is no transient indices, the stochastic matrix is said {\it completely reducible}. 

The following theorem collects some properties which are useful in the sequel.
\begin{thm}
\label{groo}
Let the stochastic matrix put in the canonical form \eqref{can}. Then
\begin{itemize}
\item[(i)] for $j=1,\dots,n$, $\s_{\rm ph}(B_{jj})$ is a subgroup of the circle group $\bt$, and the multiplicity of all peripheral eigenvalues is always 1; 
\item[(ii)] $\s(S)=\bigcup_{j=0}^n\s(B_{jj})$;
\item[(iii)] $\s_{\rm ph}(S)=\bigcup_{j=1}^n\s_{\rm ph}(B_{jj})$.
\end{itemize}
\end{thm} 
\begin{proof}
(i) follows by Theorem I.6.5 in \cite{S}, where $\s_{\rm ph}(B_{jj})$ coincides with the $d_j$-th roots of the unity, $d_j$ being the index of imprimitivity of 
$B_{jj}$, see \cite{S}, Section I.9.\footnote{If the index of imprimitivity $d_j$ of a block $B_{jj}$ in the reduced form \eqref{can} is 1, then $B_{jj}$ is said to be {\it primitive}.}

(ii) is well known, see e.g. \cite{V}, Section 2.3.\footnote{As witnessed by the matrix 
$\begin{bmatrix} 
	 1/2 & 1/4 &1/4\\
	0& 2/3&1/3\\
	0&0&1\\
\end{bmatrix}$, $B_{00}=\begin{bmatrix} 
	 1/2 & 1/4\\
	0& 2/3\\
\end{bmatrix}=\begin{bmatrix} 
	 B^{00}_{00} & B^{01}_{00}\\
	0& B^{11}_{00}\\
	\end{bmatrix}$ can be further reduced. Yet, $\s(B_{00})=\{1/3,2/3\}=\s(B^{00}_{00})\bigcup\s(B^{11}_{00})$.}

(iii) follows from (ii) because $\s_{\rm ph}(B_{00})=\varnothing$, otherwise the statement of Proposition I.9.3 in \cite{S} does not hold if $B_{00}$ has an eigenvalue 
$\l$ with $|\l|=1$.
\end{proof}

Given a stochastic matrix $S\in\bm_n(\br)$, the associated Markov chain is nothing else than a commutative finite dimensional $C^*$ dynamical system 
$(\bc^n,S)$, where the matrix $S$ generates, through its non negative powers, the action of the monoid $\bn$.

By Theorem \ref{cief}, with $\ga:=\bc^n$, $\ga_S$ which, as a linear space, coincides with $P_S\bc^n$, is in fact an abelian $C^*$-algebra with the new Choi-Effros product $\circ$. We now show the main result of the present section, that is $(\ga_S, S\lceil_{\ga_S})$ provides a genuine conservative $C^*$-dynamical system, where now the action 
of $S\lceil_{\ga_S}$ can be extended to negative powers.

Here, there is a crucial preliminary
\begin{lem}
\label{lepo}
With the above notation, $S\lceil_{\ga_S}$ is an order-automorphism.
\end{lem}
\begin{proof}
As in the proof of Proposition \ref{prtreu}, we get $S\lceil_{\ga_S}=\sum_{\l\in\s_{\rm ph}(S)}\l E_\l$. Since by Theorem \ref{groo}, 
$\s_{\rm ph}(S)=\bigcup_{j=1}^n\s_{\rm ph}(B_{jj})$ and $\s_{\rm ph}(B_{jj})$ is the cyclic subgroup of the circle consisting of the $d_j$-th roots of the unity,
$\big(S\lceil_{\ga_S}\big)^{{\rm l.c.m}\{d_1,\dots,d_n\}}=\id_{\ga_S}$. Therefore, 
$\big(S\lceil_{\ga_S}\big)^{-1}=\big(S\lceil_{\ga_S}\big)^{{\rm l.c.m}\{d_1,\dots,d_n\}-1}$ which is manifestly positive.\footnote{If all blocks $B_{jj}$, $j=1,\dots,n$, are primitive, that is $d_j=1$, then on one hand $S\lceil_{\ga_S}=\id_{\ga_S}$, and on the other hand ${\rm l.c.m}\{d_1,\dots,d_n\}-1=0$ which means
$\big(S\lceil_{\ga_S}\big)^{-1}=\id_{\ga_S}=\big(S\lceil_{\ga_S}\big)^{{\rm l.c.m}\{d_1,\dots,d_n\}-1}$.}
\end{proof}
\begin{thm}
With the above notation, $S\lceil_{\ga_S}$ is a $*$-automorphism of $\ga_S$.
\end{thm}
\begin{proof}
We need only to check that, for $a,b\in\ga_S$, $S(a\circ b)=S(a)\circ S(b)$.

By Lemma \ref{lepo}, $S\lceil_{\ga_S}$ is a order-automorphism and thus, by Theorem \ref{crstdh}, it provides a Jordan isomorphism when restricted to its selfadjoint part. Since $\ga_S$ is abelian, we have for selfadjoint elements,
$$
a\bullet b=\frac12(a\circ b+b\circ a)=a\circ b\,,
$$
and thus
$$
S(a\circ b)=S(a\bullet b)=S(a)\bullet S(b)=S(a)\circ S(b)\,.
$$

We have then recognised that $S\lceil_{\ga_S}$ preserves the product when restricted to its selfadjoint part. But $x,y\in\ga_S$ can be written as combinations of two selfadjoint elements and, since the algebra under consideration is abelian, $xy$ is written as combination of four selfadjoint elements. By linearity, $S\lceil_{\ga_S}$ preserves the product on the whole 
$\ga_S$ and the assertion follows.
\end{proof}
\begin{rem}
We want to note the following simple fact. Since $S\lceil_{\ga_S}$ is invertible, we get a, indeed reversible, $C^*$-dynamical system 
$\big(\ga_S, S\lceil_{\ga_S},\bz\big)$, with $\ga_S$ equipped with the Choi-Effros product.
\end{rem}
We end the present section with the following considerations. By taking into account Proposition \ref{vig}, point 2 in Proposition 3.2 of \cite{G}, and lastly (iii) in Theorem \ref{groo}, we conclude that in the irreducible cases, hence in all completely reducible ones, the Choi-Effros product coincides with the original one. On the other hand, we know that there are examples, necessarily with transient indices, for which the original product must be changed, see e.g. \cite{FOR}, Section 6. Therefore, one might conclude that the cases for which the original product should be changed is connected with the presence of transient indices.

Unfortunately, also this conjecture does not hold in some cases, for example when all imprimitivity indices $d_j$, $j=1,2,\dots,n$ of the block-matrices  
in the canonical form of the stochastic matrix $S$ \eqref{can} are 1, because in such a situation $\s_{\rm ph}(S)=\{1\}$ with multiplicity $n$.

Therefore, the cases for which the original product might be replaced with the Choi-Effros one are to be found among those for which the set of transient indices is nonvoid, and there exists at least one ergodic imprimitive component $B_{j_oj_o}$.

\section{some simple noncommutative examples}

The first example we want to briefly discuss is the UCP map on the full matrix algebra $\bm_n(\bc)$ obtained by pull-over of a stochastic matrix 
$S\in\bm_n(\br)$ through the conditional expectation $E_n:\bm_n(\bc)\to D_n\subset\bm_n(\bc)$, $D_n$ being the Maximal Abelian Sub-Algebra consisting of the diagonal elements. In this case, $\F=SE_n$, and $P_\F=P_SE_n$.
Therefore, $\s(\F)=\s(S)\bigcup\{0\}$, where $0$ appears with multiplicity $n(n-1)+m$, $m$ being its multiplicity in $S$.

It is easy to recognise that $\bm_n(\bc)_\F\sim \bc^n_S$, and the conservative dynamical system $\big(\bm_n(\bc)_\F, \F\lceil_{\bm_n(\bc)_\F}\big)$ associated to $\F$ is $*$-isomorphic to that $\big(\bc^n_S, S\lceil_{\bc^n_S}\big)$ both equipped, possibly, with the modified Choi-Effros product.

We now present a simple example concerning a class of UCP maps on $\bm_2(\bc)$. For the structure of completely positive maps between matrix algebras, see \cite{C}.

Indeed, let $\Phi\colon \bm_2(\mathbb{C})\to \bm_2(\mathbb{C})$ be a UCP map of the form
\[
\begin{cases}
\Phi(\,\cdot\,)= V_1^*(\,\cdot\, )V_1+V_2^*(\,\cdot\, )V_2\,,\\
V_1^*V_1+V_2^*V_2=I\,,
\end{cases}
\]
where the partial isometries $V_1,V_2$ have the form
\[
\begin{cases}
V_1(\cdot):=\langle\,\cdot\,|e_1\rangle x\,,\\
V_2(\cdot):=\langle\,\cdot\,|e_2\rangle y\,.
\end{cases}
\]
Here, $\{e_1,e_2\}$ is the canonical basis of $\mathbb{C}^2$ and $x,y\in\mathbb{C}^2$ arbitrary vectors with norm $1$ which can be written as
\[
x=e^{i\vartheta}\cos\alpha\,\,e_1+e^{i\vartheta}\sin\alpha\,\,e_2\,,\quad
y=e^{i\gamma}\cos\beta\,\,e_1+e^{i\gamma}\sin\beta\,\,e_2\,.
\]

It follows that
\begin{align*}
&\Phi\left(\begin{bmatrix}
a_{11}&a_{12}\\
a_{21}&a_{22}
\end{bmatrix}\right)\\
=&\begin{bmatrix}
a_{11}\cos^2\alpha+\frac{a_{12}+a_{21}}2\sin2\alpha+a_{22}\sin^2\alpha&0\\
0&a_{11}\cos^2\beta+\frac{a_{12}+a_{21}}2\sin2\beta+a_{22}\sin^2\beta
\end{bmatrix}\,,
\end{align*}
and thus the matrix $M_\F$ of $\Phi$ w.r.t. the basis $\{e_{11},e_{22},e_{12},e_{21}\}$ of $\bm_2(\mathbb{C})$
is
$$
M_\F=\begin{bmatrix}
\cos^2\alpha&\sin^2\alpha&\frac12\sin2\alpha&\frac12\sin2\alpha\\
\cos^2\b&\sin^2\b&\frac12\sin2\b&\frac12\sin2\b\\
0&0&0&0\\
0&0&0&0\\
\end{bmatrix}=:\begin{bmatrix}
S&B\\
O&O\\
\end{bmatrix}\,.
$$
It follows that $\s(M_\F)=\s(S)\bigcup\{0\}$ and $\s_{\rm ph}(M_\F)=\s_{\rm ph}(S)$.

Suppose now that $\l\in\s(M_\F)\smallsetminus\{0\}$. With $\xi=\begin{bmatrix}
a_{11}\\
a_{22}\\
\end{bmatrix}$ and $\eta=\begin{bmatrix}
a_{12}\\
a_{21}\\
\end{bmatrix}$, $M_\F\begin{bmatrix}
\xi\\
\eta\\
\end{bmatrix}=\l \begin{bmatrix}
\xi\\
\eta\\
\end{bmatrix}$ reads
$$
S\xi+B\eta=\l\xi\,\,\,\&\,\,\, \l\eta=0\,,
$$
which leads to
$$
\eta=0 \,\,\,\&\,\,\, S\xi=\l\xi\,.
$$

Therefore, even if the situation of the latter example is different from the former one, at the same time we have $\bm_2(\bc)_\F\sim \bc^2_S$, with 
$\big(\bm_2(\bc)_\F, \F\lceil_{\bm_2(\bc)_\F}\big)$ being $*$-isomorphic to $\big(\bc^2_S, S\lceil_{\bc^2_S}\big)$ both equipped, possibly, with the modified Choi-Effros product.

\end{document}